\documentclass[letterpaper, 10 pt, conference]{ieeeconf}  % Comment this line out
\IEEEoverridecommandlockouts                              % This command is only
\usepackage{cite, graphicx, amsmath, amssymb, mathrsfs, float, subfig, color, url, tcolorbox, mathpazo}
\usepackage{etoolbox}

\newtheorem{theorem}{{\bf Theorem}}

\renewcommand{\baselinestretch}{1}

\title{\LARGE \bf
Robust Distributed Control of DC Microgrids with Time-Varying Power Sharing
}

\author{Mayank Baranwal$^{1,a}$, Alireza Askarian$^{1,b}$, Srinivasa M. Salapaka$^{1,c}$ and Murti V. Salapaka$^{2,d}$% <-this % stops a space
\thanks{$^{1}$Department of Mechanical Science and Engineering, University of Illinois at Urbana-Champaign, 61801 IL, USA}
\thanks{$^{1}$Department of Electrical and Computer Engineering, University of Minnesota, Minneapolis, 55455 MN, USA}
\thanks{$^{a}${\tt\small baranwa2@illinois.edu}, $^{b}${\tt\small askaria2@illinois.edu}, $^{c}${\tt\small salapaka@illinois.edu}, $^{d}${\tt\small murtis@umn.edu}}
}

\begin{document}

\maketitle
\thispagestyle{empty}
\pagestyle{empty}

\begin{abstract}
This paper addresses the problem of output voltage regulation for multiple DC/DC converters connected to a microgrid, and prescribes a scheme for sharing power among different sources. This architecture is structured in such a way that it admits quantifiable analysis of the closed-loop performance of the network of converters; the analysis simplifies to studying closed-loop performance of  an equivalent {\em single-converter} system. The proposed architecture allows for the proportion in which the sources provide power to vary with time; thus overcoming limitations of our previous designs in \cite{baranwal2016robust}. Additionally, the proposed control framework is suitable to both centralized and decentralized implementations, i.e., the same control architecture can be employed for voltage regulation irrespective of the availability of common load-current (or power) measurement, without the need to modify controller parameters. The performance becomes quantifiably better with better communication of the demanded load to all the controllers at all the converters (in the centralized case); however guarantees viability when such communication is absent. Case studies comprising of battery, PV and generic sources are presented and demonstrate the enhanced performance of prescribed optimal controllers for voltage regulation and power sharing.
\end{abstract}

%%%%%%%%%%%%%%%%%%%%%%%%%%%%%%%%%%%%%%%%%%%%%%%%%%%%%%%%%%%%%%%%%%%%%%%%%%%%%%%%
\section{INTRODUCTION}\label{sec:Intro}
\begin{figure}[!t]
	\centering
	\includegraphics[width=0.95\columnwidth]{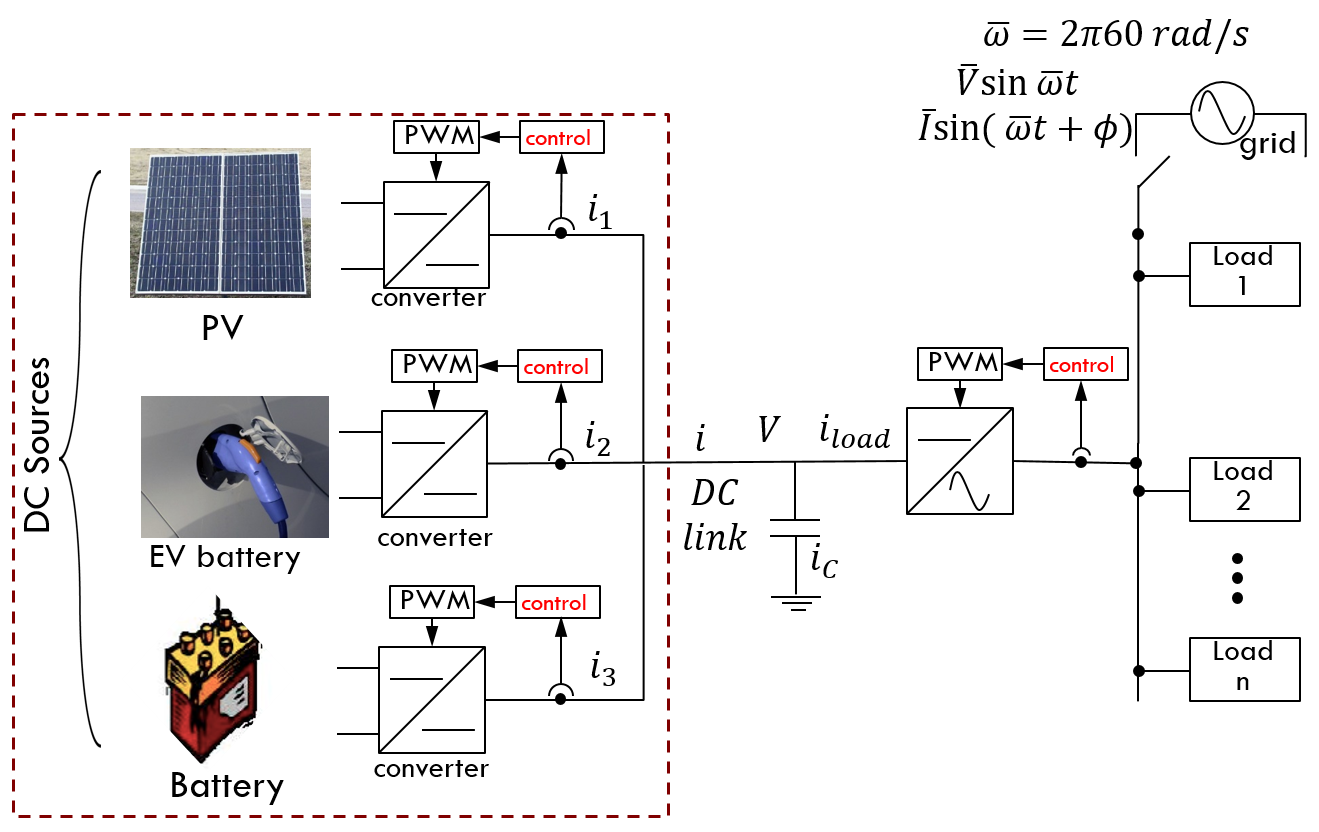}
	\vspace{-0.5em}
	\caption{{\small A schematic of a microgrid. An array of DC sources provide power for AC loads. Power sources  provide  power at DC-link, their common output bus, at a voltage that is regulated to a set-point. The control system at the respective DC-DC converter that interfaces with a source is responsible for regulating the voltage at the DC-link. An inverter that connects to the DC-link converts the total current from the sources at the regulated voltage to  alternating current (AC) at its output to satisfy the power demands of the AC loads. This paper describes an approach for control design of the multiple converters systems associated with power transfer from sources to the DC-link (shown by the dotted line).}} 
	\label{fig:MCS}
	\vspace{-2em}
\end{figure}

Utility grid is most stressed during peak power demands resulting in significant increase in real-time power prices and congestion in the local power distribution zone. Microgrids can help reduce the requirement for additional utility generation and thus minimize the demand on the utility grid by enabling integration of renewable energy sources such as solar and wind energy, distributed energy resources (DERs), energy storage, and demand response. Microgrids are localized grid systems that are capable of operating in parallel with, or independently from, the existing traditional grid \cite{lasseter2002microgrids, hatziargyriou2007microgrids}. Fig. \ref{fig:MCS} shows a schematic of a microgrid with multiple DC sources providing power for AC loads. Existing control architectures for traditional grids, which are  designed for relatively large conventional sources (power plants) of  predictable and {\em dispatchable} electric power, cannot adequately manage {\em uncertain} power sources such as solar or wind generations. Limited predictability with such resources result in intermittent power generation; moreover time-varying loads, practicability and economics factors pose additional challenges in efficient operation of microgrids. Thus it is required to develop efficient distributed control technologies for reliable operation of smart microgrids \cite{lopes2006defining}.

In such smart grids, multiple DC power sources connected in parallel, each interfaced with DC-DC converter, provide power at their common output, the DC-link, at a regulated voltage; this power can directly feed DC loads or be used by an inverter to interface with AC loads (see Fig. \ref{fig:MCS}). By appropriately controlling the switch duty-cycle of DC-DC converter at each power source, it becomes possible to manipulate electrical quantities such as the power output by the power source and the   voltage at the DC-link. The main goals of the control design is to regulate voltage at the DC-link and ensuring a prescribed sharing of power between different sources; for instance, economic considerations can dictate that power provided by the sources should be in a certain proportion or according to a prescribed priority (e.g. PV$>$battery$>$EV where PV provides the maximum power it can to satisfy load demand, and the deficit is provided by battery,and then the EV). The main challenges arise from the uncertainties in the size and the schedules of loads, the complexity of a coupled multi-converter network, the uncertainties in the model parameters at each converter, and the adverse effects of interfacing DC power sources with AC loads, such as the $120$ Hz ripple that has to be provided by the DC sources.

Problems pertaining to robust and optimal control of converters have received recent attention. Conventional PID-based controllers often fail to address the problem of robustness and modeling uncertainties. In \cite{olalla2010lmi}, a linear-matrix-inequality (LMI) based robust control design is presented for boost converters which demonstrates significant improvements in voltage regulation over PID based control designs. In \cite{weiss2004h, hornik2011current, salapaka2014viability}, robust $\mathcal{H}_\infty$-control framework is employed in the context of inverter systems. While the issue of current sharing is extensively studied \cite{panov1997analysis, wu1993load}, most methods assume a single power source. A systematic control design that addresses all the challenges and objectives for the multi-converter control is still lacking.

The control architecture proposed in this paper addresses the following primary objectives - a) voltage regulation at the DC-link with guaranteed robustness margins, b) prescribed {\em time-varying} power sharing in a network of parallel converters, c) controlling the trade-off between $120$Hz ripple on the total current provided by the power sources and the ripple on the DC-link voltage. While these objectives are partially addressed in our prior work \cite{baranwal2016robust} on the robust control of DC-DC converters, a main drawback of the design proposed in \cite{baranwal2016robust} is that the control framework does not allow for time-varying power sharing requirements. In this work, we propose a new architecture wherein the power requirements on each converter are imposed through external reference signals; this allows for time-varying power sharing/priority prescriptions. This is achieved without sacrificing any advantages of the design in \cite{baranwal2016robust}.

An important feature of the proposed architecture is that it exploits structural features of the paralleled multi-converter system, which results in a modular and yet coordinated control design. Accordingly at each converter, it employs a nested (outer-voltage inner-current) control structure \cite{erickson2007fundamentals}, where all converters share the same design for the outer-loop voltage controllers while the inner-loop current controllers are so chosen that the entire closed-loop multi-converter system can be reduced to an equivalent {\em single-converter} system in terms of the transfer function from the desired regulation setpoint $V_\text{ref}$ to the voltage at the DC-link $V$. An interesting aspect of the proposed implementation is that the same control implementation works for both the {\em centralized} case, when the  load power is known and communicated  to all the controllers, and the {\em decentralized} case where the load is unknown. The architecture achieves better performance in voltage regulation and power sharing when load power  information is communicated, while it guarantees electrical viability with quantifiable bounds on deviations from the targeted performances in the decentralized case.

\section{MODELING OF CONVERTERS}\label{sec:model}
\begin{figure}[!t]
	\centering
	\includegraphics[width=0.9\columnwidth]{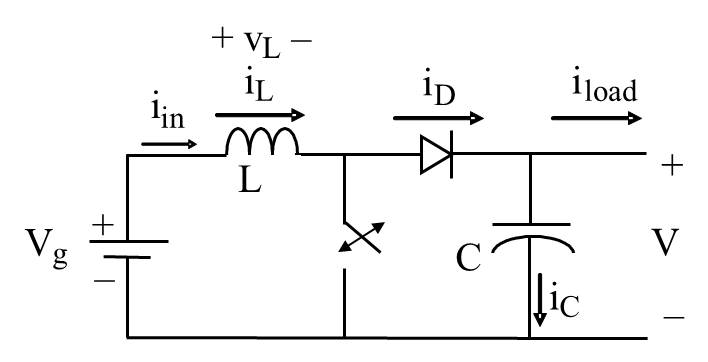}
	\vspace{-0.5em}
	\caption{{\small A schematic of a Boost converter. The converters are assumed to operate in continuous-conduction-mode (CCM). $V_g$ is the voltage of input source such as EV, battery and PV. The DC-link voltage $V$ is controlled using a switch by varying its duty-cycle. Note that the input source current $i_\text{in} = i_L$.}} 
	\label{fig:conv_schematic}
	\vspace{-2em}
\end{figure}
In this section, we describe the differential equations that govern the dynamics of DC-DC converters. These converters belong to a class of switched-mode power electronics, where a semiconductor based high-frequency switching mechanism (and associated electronic circuit) connected to a DC power source enables changing voltage and current characteristics at its output. The models presented below depict dynamics for signals that are averaged over a switch cycle.

Fig. \ref{fig:conv_schematic}a shows a schematic of a Boost converter. Boost converter regulates a voltage $V$ at its output which is larger than the input voltage $V_g$. The averaged dynamic model of a Boost converter is given by
\begin{small}
	\begin{eqnarray}\label{eq:boost_eq}
		L\dot{i}_L(t) &=& -(1-d(t))V(t) + V_g,\nonumber \\
		C\dot{V}(t) &=& (1-d(t))i_L(t) - i_\text{load}(t),
	\end{eqnarray}
\end{small}
where $d(t)$ represents the duty-cycle (or the proportion of {\em ON} duration) at time $t$. By defining $d'(t):=1-d(t)$ and $D':=\left(V_g/V_\text{ref}\right)$, where $V_\text{ref}$ is the desired output voltage, (\ref{eq:boost_eq}) can be rewritten as,
\begin{small}
	\begin{eqnarray}\label{eq:boost_eq1}
		L\frac{di_L(t)}{dt} &=& \underbrace{V_g - d'(t)V(t)}_{\tilde{u}(t):=V_g-u(t)},\nonumber \\
		C\frac{dV(t)}{dt} &=& \underbrace{(D'+\hat{d}(t))}_{\approx D'}i_L(t) - i_\text{load}(t).
	\end{eqnarray}
\end{small}
Note that $\hat{d}(t) = d'(t)-D'$ is typically very small, and therefore allows for a linear approximation around the nominal duty-cycle, $D = 1-D'$ (see \cite{erickson2007fundamentals} for details).

Similarly we can describe the averaged dynamics of a Buck converter,
\begin{small}
	\begin{eqnarray}\label{eq:buck_eq}
		L\frac{di_L(t)}{dt} &=& \underbrace{-V(t) + d(t)V_g}_{\tilde{u}(t):=-V(t)+u(t)} = \tilde{u}(t),\nonumber \\
		C\frac{dV(t)}{dt} &=& i_L(t) - i_\text{load}(t).
	\end{eqnarray}
\end{small}

\section{PROBLEM FORMULATION}\label{sec:problem}
This paper addresses the following primary objectives {\em simultaneously} (in context of Fig. \ref{fig:MCS}) - (a) Regulation of DC-link voltage $V$ to a prescribed value $V_{\text{ref}}$ in presence of time-varying loads/generation and parametric uncertainties, (b) time-varying current (power) sharing among multiple sources, that is, ensuring that current (power) outputs $i_k$ from each converter respectively track a time-varying signal $i_{\text{ref}_k}$, and (c) $120 Hz$ ripple current sharing between the output currents $i_k$ from each converter and the capacitor current $i_C$. The last objective is dealt in our prior work \cite{baranwal2016robust} and is addressed by an appropriate design of inner-controller described in Sec. \ref{subsec:Single_inner} and \ref{sec:Control_Many}. In this paper, we primarily focus on achieving the first two objectives, while inheriting the properties of the inner-controller for ripple current sharing.

\section{CONTROL FRAMEWORK FOR SINGLE CONVERTER}\label{sec:Control_Single}
\begin{figure}
	\includegraphics[width= 0.9\columnwidth]{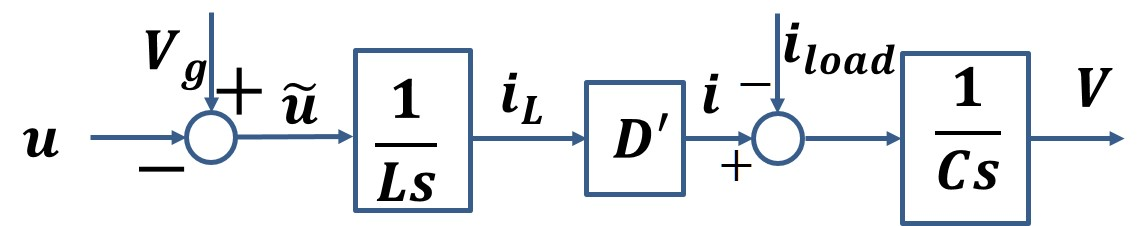}
	\caption{{\small Block diagram representation of a boost-type converter. The control signal $\tilde{u}$ is converted to an equivalent PWM signal to command the gate of the transistor acting as a switch.}}
	\vspace{-1.5em}
	\label{fig:boost_model}
\end{figure}

In this section, we describe the {\em inner-outer} control design for a single Boost converter system. This design for a {\em single} converter forms the basis for the analysis and design of control architecture for multiple converters presented in Sec. \ref{sec:Control_Many}. Fig. \ref{fig:boost_model} shows the block diagram representation of the Boost converter based on the dynamical equations in (\ref{eq:boost_eq1}). While the design is easily extendable to include other converter types such as Buck and Buck-Boost, the discussion has been confined to Boost converters only for the sake of brevity. Note that in the proposed control architecture (see Fig. \ref{fig:single_block}), the inputs to outer feedback controller include $i_{\text{ref}}$ in addition to the typical $V_{\text{ref}}$ and the measured DC-link voltage $V$. The requirements on current sharing are imposed through this additional $i_{\text{ref}}$ signal (explained in Sec. \ref{sec:Control_Many}) by setting $i_{\text{ref}}$ to measured (or communicated) load current $i_{\text{load}}$ in the centralized case, and setting $i_{\text{ref}}$ to estimated (or prespecified signals) in the decentralized case.

\subsection{Design of the inner-loop controller}\label{subsec:Single_inner}
The design for the inner-loop controller $K_c$ is inherited from our previous work \cite{baranwal2016robust}. The main objective for designing the inner-loop controller $K_c$ is to decide the trade-off between the $120$Hz ripple on the capacitor current $i_C$ (equivalently on the output voltage $V$) and the inductor current $i_L$ of the converter. Accordingly, $K_c$ is designed such that the inner-shaped plant $\tilde{G}_c$ is given by
\begin{small}
	\begin{equation}\label{eq:Gc_tilde}
		\tilde{G}_c(s) = \left(\frac{\tilde{\omega}}{s+\tilde{\omega}}\right)\left(\frac{s^2+2\zeta_1\omega_0s+\omega_0^2}{s^2+2\zeta_2\omega_0s+\omega_0^2}\right).
	\end{equation}
\end{small} 
where $\omega_0 = 2\pi 120$rad/s and $\tilde{\omega}, \zeta_1, \zeta_2$ are design parameters. The parameter $\tilde{\omega} > \omega_0$ and it is used to implement a low-pass filter to attenuate undesirable frequency content in $i_L$ beyond $\tilde{\omega}$. Thus, the bandwidth of the inner closed-loop plant is decided by the choice of $\tilde{\omega}$. The parameters $\zeta_1$ and $\zeta_2$ impart a notch-like behavior to $\tilde{G}_c$ at $\omega_0 = 120$Hz, and the size of the notch is determined by the ratio $\zeta_1/\zeta_2$. Note that $\tilde{G}_c$ represents the inner closed-loop plant from the output of the outer-loop controllers $\tilde{u}$ to the inductor current $i_L$, and since $i_C = i_L - i_\text{load}$, the ratio $\zeta_1/\zeta_2$ can be appropriately designed to achieve a specified trade-off between $120$Hz ripple on $i_C$ and $i_L$, and therefore between $V$ and $i_L$. The stabilizing $2^{nd}$-order controller $K_c$ that yields the aforementioned inner closed-loop plant $\tilde{G}_c$ is explicitly given by
\begin{small}
	\begin{equation}\label{eq:Kc}
		K_c(s) = L\tilde{\omega}\frac{(s^2+2\zeta_1\omega_0s+\omega_0^2)}{(s^2+2\zeta_2\omega_0s+\omega_0^2+2(\zeta_2-\zeta_1)\omega_0\tilde{\omega})}.
	\end{equation}
\end{small}
The readers are encouraged to refer to Sec. III in \cite{baranwal2016robust} for further details on the inner-loop control design.

\subsection{Design of the outer-loop controller}\label{subsec:Single_outer}
\begin{figure}
	\includegraphics[width= 0.95\columnwidth]{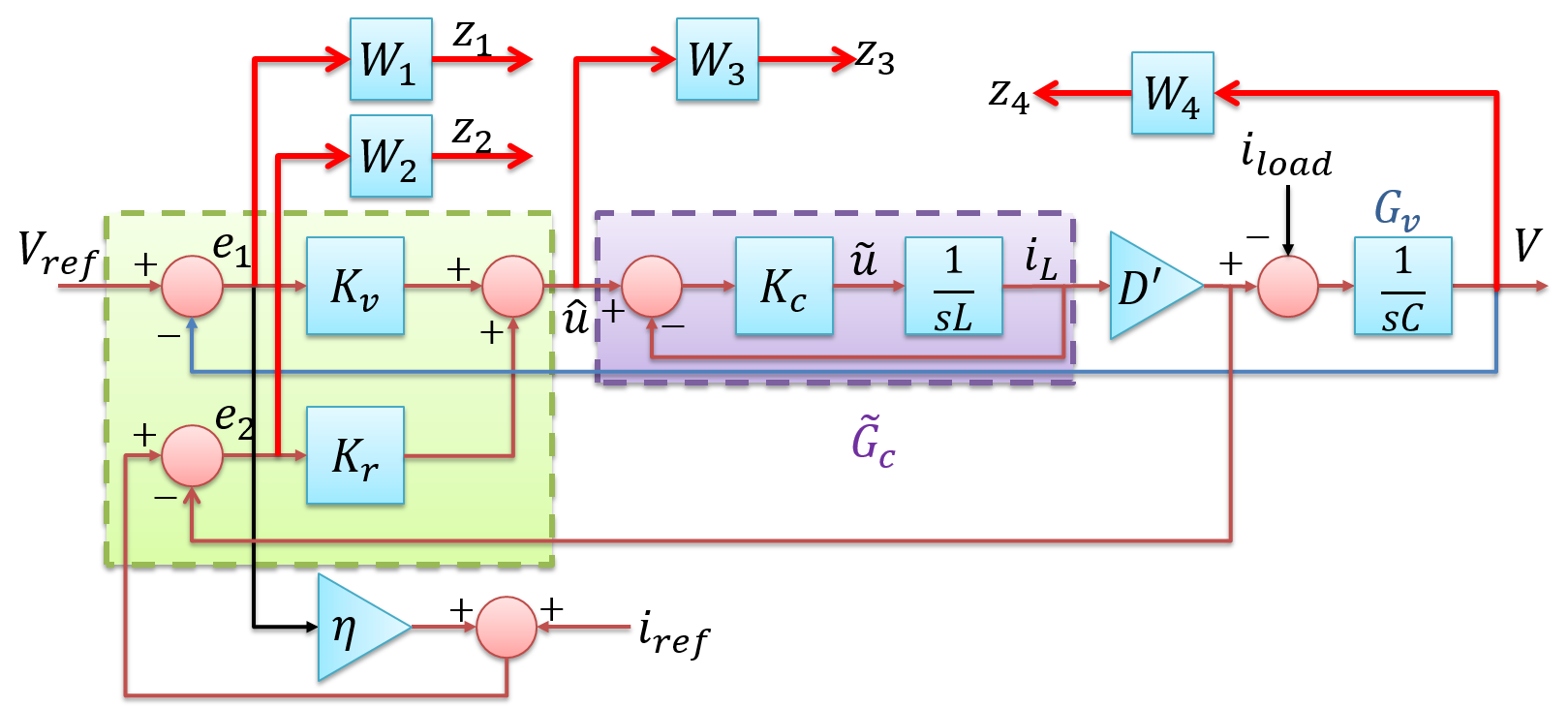}
	\caption{{\small Block diagram representation of the inner-outer control design. Exogenous signal $V_\text{ref}$ represents the desired output voltage. The quantities $V$, $i_\text{load}$ and $i_L$ represent the output voltage, load current and inductor current, respectively. The {\em regulated variables} $z_1, z_2, z_3$ and $z_4$ correspond to weighted - (a) tracking error in DC-link voltage, (b) mismatch between $i_\text{ref}$ and $i_\text{load}$, (c) control effort $\hat{u}$, and (d) output voltage tracking, respectively.}}
	\vspace{-1.5em}
	\label{fig:single_block}
\end{figure}
For a given choice of inner-controller $K_c$, we present our analysis and design of controller in terms of transfer function block diagrams shown in Fig. \ref{fig:single_block}. In this figure, $\tilde{G}_c$ represents the inner shaped plant. The outer controllers are denoted by $K_{v}$ and $K_{r}$, and are designed to regulate the output DC voltage $V$ to the desired reference voltage $V_\text{ref}$ and the output current $D'i_L$ to the reference current $i_\text{ref}$, respectively. Note that from (\ref{eq:boost_eq1}), $D'i_L$ is equal to $i_\text{load}$ at steady-state. The augmentation of controller $K_{r}$ forms the basis for time-varying power sharing and is explained in the next section. It should be remarked that the proposed design has a feature that if the load current measurement is available, i.e., $i_\text{ref} = i_\text{load}$, then the steady state DC output voltage is maintained at the reference voltage $V_\text{ref}$. However in the absence of $i_\text{load}$ measurement, the outer controller $K_{r}$ regulates the output current $D'i_L$ to $i_\text{ref}\neq i_\text{load}$ resulting in an output voltage $V\neq V_\text{ref}$. The mismatch in voltage tracking is captured by a pre-specified {\em droop}-like coefficient $\eta$ in a controlled manner, the notable difference here being the application of droop to the {\em faster} current loop when compared with the conventional droop-based design acting on the {\em slower} voltage loop. This feature is mathematically quantified in the following discussion on the proposed control design.

The performance of a Boost converter is characterized by its voltage and power reference tracking bandwidths, better voltage signal-to-noise-ratio (SNR), and robustness to modeling uncertainties. The main objective for the design of the controllers $K_{v}$ and $K_{r}$ is to make the tracking errors small and simultaneously attenuate measurement noise to achieve high resolution. This is achieved by posing a {\em model-based} multi-objective optimization problem, where the required objectives are described in terms of norms of the corresponding transfer functions, as described below. Note that the regulated variables $z_1, z_2, z_3$ and $z_4$ correspond to weighted - (a) tracking error in DC-link voltage, (b) mismatch between $i_\text{ref}$ and $i_\text{load}$, (c) control effort $\hat{u}$, and (d) output voltage tracking, respectively. From Fig. \ref{fig:single_block}, the transfer function from exogenous inputs and auxiliary control input $w = \left[V_\text{ref}, i_\text{ref}, i_\text{load}, \hat{u}\right]^T$ to regulated output $z = \left[z_1, z_2, z_3, z_4, e_1, e_2\right]$ is given by
\begin{small}
	\begin{equation}\label{eq:hinf}
	\resizebox{.5 \textwidth}{!} 
{
	$\underbrace{\left[ \begin{array}{c}
		z_1\\
		z_2\\
		z_3\\
		z_4\\
		e_1\\
		e_2
	\end{array} \right]}_{z} = 
	\underbrace{\left[ \begin{array}{cccc}
		W_1 & 0 & W_1G_v & -D'W_1G_v\tilde{G}_c \\
		\eta W_2 & W_2 & \eta W_2G_v & -D'(1+\eta G_v)W_2\tilde{G}_c\\
		0 & 0 & 0 & W_3 \\
		0 & 0 & -W_4G_v & D'W_4G_v\tilde{G}_c \\
		1 & 0 & G_v & -D'G_v\tilde{G}_c \\
		\eta & 1 & \eta G_v & -D'(1+\eta G_v)\tilde{G}_c
	\end{array}	\right]}_{T_{wz}} 
	\underbrace{\left[ \begin{array}{c}
		V_\text{ref}\\
		i_\text{ref}\\
		i_\text{load}\\
		\hat{u}
	\end{array} \right]}_w.$
}
	\end{equation}
	
\end{small}
The optimization problem is to find stabilizing controllers $K_{outer} = \left[K_{v},K_{r}\right]^T$ such that the $\mathcal{H}_\infty$-norm of the above transfer function from $w$ to $z$ is minimized. Here the weights $W_1, W_2, W_3$ and $W_4$ are chosen to reflect the design specifications of robustness to parametric uncertainties, tracking bandwidth, and saturation limits on the control signal. More specifically, the weight functions $W_1(j\omega)$ and $W_2(j\omega)$ are chosen to be large in frequency range $\left[0,\omega_{BW}\right]$ to ensure small tracking errors $e_1 = V_\text{ref}-V$ and $e_2 = i_\text{ref} + \eta e_1 - D'i_L$ in this frequency range. The design of weight function $W_3(j\omega)$ entails ensuring that the control effort lies within saturation limits. The weight function $W_4$ is designed as a high-pass filter to ensure that the transfer function from $i_\text{load}$ to $V$ is small at high frequencies to provide mitigation to measurement noise.

Note that for the system shown in Fig. \ref{fig:single_block}, the voltage $V$ at the DC-link is given by,
\begin{small}
	\begin{equation}\label{eq:V_eq}
		V = G_v\left(-i_\text{load} + D'\tilde{G}_c(K_ve_1 + K_re_2)\right).
	\end{equation}
\end{small}
Using the fact that $e_1=V_\text{ref}-V$ and $e_2=i_\text{ref}+\eta e_1-D'\tilde{G}_c(K_ve_1 + K_re_2)$, the DC-link voltage in terms of exogenous quantities $V_\text{ref}, i_\text{ref}$ and $i_\text{load}$ is given by
\begin{small}
	\begin{eqnarray}\label{eq:voltage_single}
		&&V(s) = \left[\frac{D'\tilde{G}_cG_v(K_v+\eta K_r)}{1+D'\tilde{G}_cK_r+D'\tilde{G}_cG_v(K_v+\eta K_r)}\right]V_\text{ref}(s) \nonumber\\
		&&+ \left[\frac{D'\tilde{G}_cG_vK_r}{1+D'\tilde{G}_cK_r+D'\tilde{G}_cG_v(K_v+\eta K_r)}\right]\left(i_\text{ref}(s)-i_\text{load}(s)\right) \nonumber\\
		&&- \left[\frac{G_v}{1+D'\tilde{G}_cK_r+D'\tilde{G}_cG_v(K_v+\eta K_r)}\right]i_\text{load}(s).
	\end{eqnarray}
\end{small}
Let $S(s), T_{V_\text{ref}V}$ and $T_{i_\text{ref}V}$ denote the closed-loop sensitivity transfer function and complementary sensitivity transfer functions from $V_\text{ref}$ to $V$ and $i_\text{ref}$ to $V$, respectively. Then (\ref{eq:voltage_single}) can be rewritten as
\begin{small}
	\begin{eqnarray}
		V(s) = T_{V_\text{ref}V}V_\text{ref}(s) + T_{i_\text{ref}V}\left(i_\text{ref}(s)-i_\text{load}(s)\right) - G_vSi_\text{load}(s)\nonumber.
	\end{eqnarray}
\end{small}
The DC gains of above closed-loop transfer functions are given by (since $G_v =1/sC$ has an infinite DC gain),
\begin{small}
	\begin{eqnarray}
		&&|T_{V_\text{ref}V}(j0)|=1,\quad |T_{i_\text{ref}V}(j0)|=\frac{|K_r(j0)|}{|K_v(j0)+\eta K_r(j0)|}\quad \text{and}\nonumber\\
		&&|(G_vS)(j0)| = \frac{1}{D'\left(|K_v(j0)+\eta K_r(j0)|\right)}.\nonumber
	\end{eqnarray}
\end{small}
We now provide a sketch of the proposed design concept. Since $K_v$ and $K_r$ are chosen as high DC-gain controllers (obtained by solving the $\mathcal{H}_\infty$ optimization problem), we have $|G_vS((j\omega))|\approx 0$ at low-frequencies. Thus the effect of disturbance signal $i_\text{load}$ is insignificant at low frequencies. Similarly $T_{V_\text{ref}V}(j\omega)$ has unity gain at low frequencies. Furthermore, if the load current $i_\text{load}$ measurement is available (i.e. $i_\text{ref} = i_\text{load}$), then the Boost converter tracks the reference voltage with almost unity gain. However in the absence of $i_\text{load}$ measurement, the tracking error depends on the mismatch between $i_\text{ref}$ and $i_\text{load}$, i.e., the bound on the steady-state tracking error becomes proportional to $\dfrac{|K_r(j0)|}{|K_v(j0)+\eta K_r(j0)|}$ multiplied by the mismatch value $|i_\text{ref}(j0)-i_\text{load}(j0)|$. By choosing appropriate controllers $K_v$ and $K_r$ (i.e., $|T_{i_\text{ref}V}(j0)|<<1$), the tracking error can be made small.

{\em Extension to Buck Converters}: The extension of the proposed control design to Buck and Buck-Boost DC-DC converters is easily explained after noting that their averaged models are structurally identical to Boost converters, except that the dependence of duty cycles on the control signal $u$ or constant parameter $D'$ are different. The differences in how duty cycles depend on $u(t)$ do not matter from the control design viewpoint since duty cycles for pulse-width modulation are obtained only after obtaining the control designs (that use the averaged models).

\section{EXTENSION TO A SYSTEM OF PARALLEL CONVERTERS}\label{sec:Control_Many}
\begin{figure*}[!t]
	\begin{center}
	\begin{tabular}{cc}
	\includegraphics[width=0.86\columnwidth]{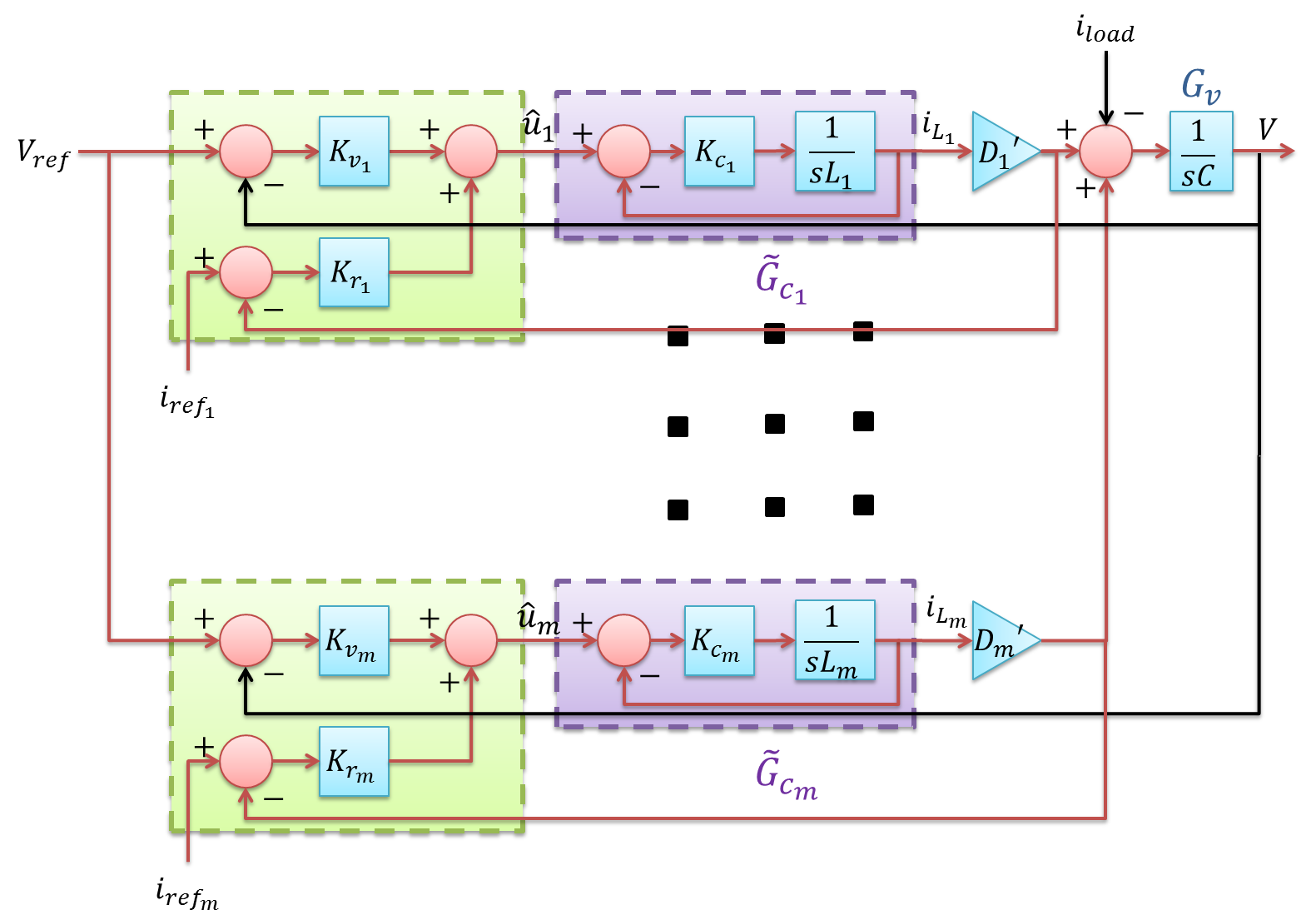}&\includegraphics[width=0.95\columnwidth]{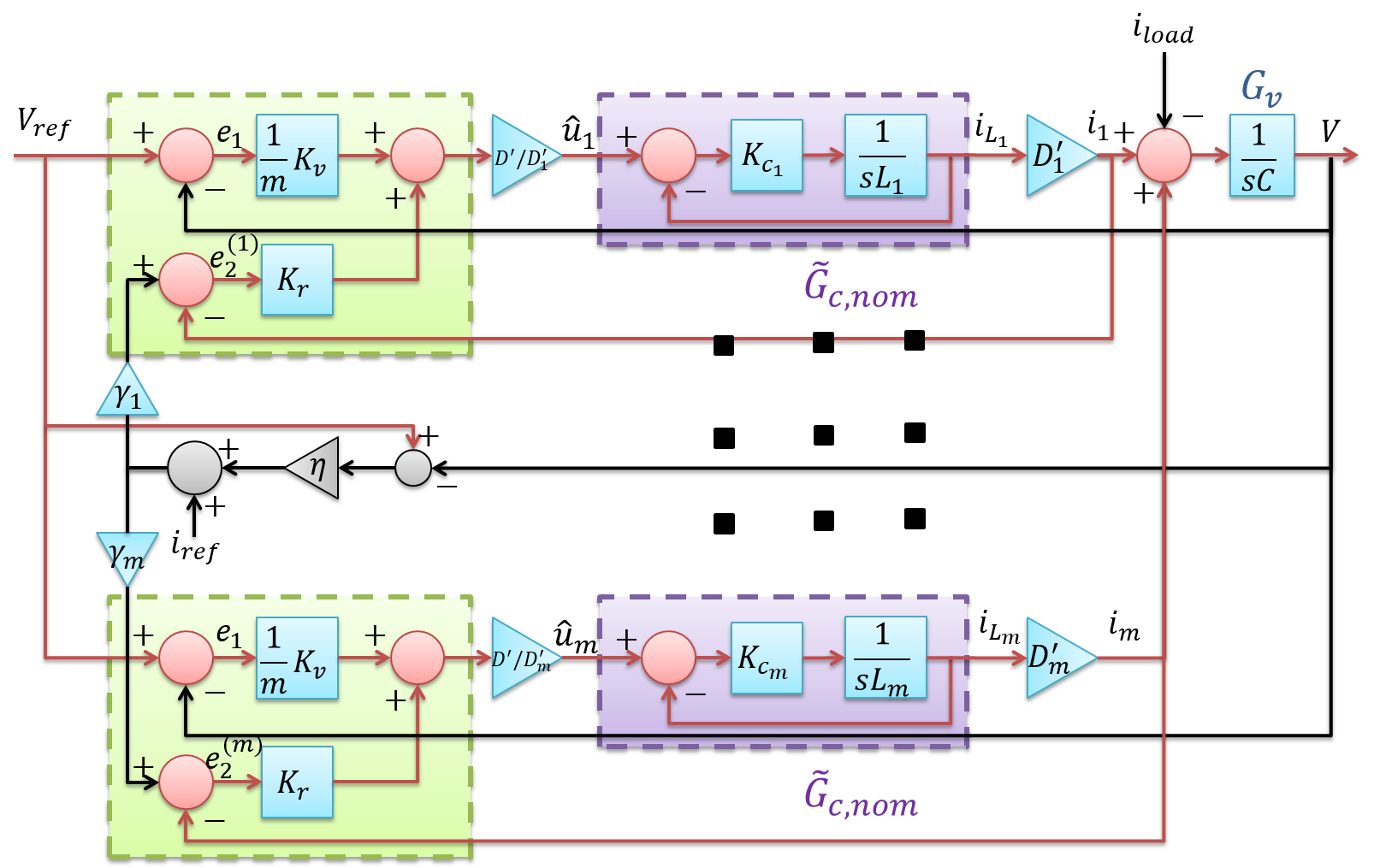}\cr
	(a)&(b)\end{tabular}
	\vspace{-0.5em}
	\caption{{\small (a) General control framework for a network of $m$ parallel converters. (b) A multiple-converters system with shaped inner plants $\tilde{G}_{c}$. In the proposed implementation, we adopt the same outer controller for different converters, i.e., $K_{v_1} = K_{v_2} = .. = K_{v_m} = \frac{1}{m}K_v$ and $K_{r_1} = K_{r_2} = .. = K_{r_m} = K_r$. Here $\gamma_k$ represents the proportion of power demanded from the $k^{th}$ source.}}
	\label{fig:model_many}
	\vspace{-2em}
	\end{center}
\end{figure*}

In this section	we extend our control framework for a single converter to a system of DC-DC converters connected in parallel in the context of power sharing, keeping in mind the practicability and robustness to modeling and load uncertainties. In particular, we analyze the multi-converter system in Fig. \ref{fig:model_many}b through an {\em equivalent} single-converter system (similar to the system shown in Fig. \ref{fig:single_block}), where the multi-converter system inherits the performance and robustness achieved by a design for the single-converter system.

Fig. \ref{fig:model_many}a represents an inner-outer control framework for a system of $m$ parallel connected converters. Note that instead of feeding the reference current $i_\text{ref$_k$}$ directly to the $k^{th}$ outer controller $K_{r_k}$, the reference signal $i_\text{ref}+\eta(V_\text{ref}-V)$ is prescaled by a time-varying multiplier $\gamma_k, 0\leq\gamma_k\leq 1$. The choice of $\gamma_k$ dictates the power sharing requirements on the $k^{th}$ converter. In fact, we later show that the proposed implementation distributes the output power in the ratios $\gamma_1:\gamma_2:..:\gamma_m$. After noting that the voltage-regulation and current reference tracking is common to all the outer controllers, in our architecture, we impose the same design for outer-controllers for all the converters, i.e., $K_{v_1} = K_{v_2} = .. = K_{v_m}$ and $K_{r_1} = K_{r_2} = .. = K_{r_m}$. This imposition enables significant reduction in the overall complexity of the distributed control design for a parallel network of converters and power sources, thus ensuring the practicability of the proposed design which allows integration of power sources of different types and values.

We design inner-controllers $K_{c_k}$ such that the inner-shaped plants from $\tilde{u}_k$ to $i_{L_k}$ are same and given by,
\begin{small}
	\begin{equation}\label{eq:Gcn}
		\tilde{G}_{c,nom}(s) = \left(\frac{\tilde{\omega}}{s+\tilde{\omega}}\right)\left(\frac{s^2+2\zeta_{1,nom}\omega_0s+\omega_0^2}{s^2+2\zeta_{2,nom}\omega_0s+\omega_0^2}\right),
	\end{equation}
\end{small}
where the ratio $\zeta_{1,nom}/\zeta_{2,nom}$ determines the tradeoff of $120$Hz ripple between the total output current $D'i_L = \sum\limits_{k=1}^mD'_ki_{L_k}$ and the capacitor current $i_C$. Note that for given values of $\zeta_{1,nom}, \zeta_{2,nom}$ and inductance $L_k$, explicit design of $K_{c_k}$ exists and is given by (\ref{eq:Kc}). After noting that $K_{v_k} = \frac{1}{m}K_{v}$ and $K_{r_k} = K_{r}$, the system in Fig. \ref{fig:model_many}a can be simplified to Fig. \ref{fig:model_many}b.

Indeed, by our choice of inner and outer controllers, the transfer functions from external references $V_\text{ref}$, $i_\text{ref}$ and $i_\text{load}$ to the desired output $V$ are identical for all converters. Hence the entire network of parallel converters can be analyzed in the context of an equivalent single converter system. This implies that $K_{v_k}$ and $K_{r_k}$ can be computed by solving $\mathcal{H}_\infty$-optimization problem (as discussed in the previous section) similar to the {\em single} converter case. We make these design specifications more precise and bring out the equivalence of the control design for the single and multiple converter systems in the following theorem.

We say that the system representation in Fig. \ref{fig:single_block} is {\em equivalent} to that in Fig. \ref{fig:model_many}b, when the transfer functions from the reference voltage $V_\text{ref}$, reference current $i_\text{ref}$ and load current $i_\text{load}$ to the DC-link voltage $V$ in Fig. \ref{fig:single_block} are identical to the corresponding transfer functions in Fig. \ref{fig:model_many}b. 

\begin{theorem}\label{thm1}
Consider the single-converter system in Fig. \ref{fig:single_block} with inner-shaped plant $\tilde{G}_{c,nom}(s)$ as given in (\ref{eq:Gcn}), outer controllers $K_{v}$, $K_{r}$, droop-coefficient $\eta$, and external references $V_\text{ref}, i_\text{load}, i_\text{ref}$; and the multi-converter system described in Figs. \ref{fig:model_many}a and \ref{fig:model_many}b with inner-shaped plants $\tilde{G}_{c_k}=\tilde{G}_{c,nom}(s)$ and outer controllers $K_{v_k}=\frac{1}{m}K_v$; $K_{r_k}=K_{r}$, droop-coefficient $\eta$, and same external references $V_\text{ref}, i_\text{load}$ and reference current $i_\text{ref}$ prescaled by time-varying scalars $\gamma_k>0$ for $1\leq k\leq m$.\\
\noindent {\bf 1. [System Equivalence]}: If $\sum_{k=1}^m\gamma_k = 1$, then the system representation in Fig. \ref{fig:single_block} is {\em equivalent} to the system representation in Fig. \ref{fig:model_many}b.\\
\noindent {\bf 2. [Power Sharing]}: For any two converters $k$ and $l$, $k,l\in\{1,\hdots,m\}$ in a multi-converter system shown in Fig. \ref{fig:model_many}b, the difference in the corresponding steady-state scaled output currents is given by
\begin{equation}\label{eq:thm1}
	\resizebox{.48 \textwidth}{!} 
	{
		$\left|\frac{D_k'i_{L_k}(j0)}{\gamma_k}-\frac{D_l'i_{L_l}(j0)}{\gamma_l}\right|\leq \left(\eta|\tilde{T}_1(j0)| + \left|\frac{1}{\gamma_k}-\frac{1}{\gamma_l}\right||\tilde{T}_2(j0)|\right)|e_1(j0)|$,
	}
\end{equation}
where, $\tilde{T}_1 := \frac{D'\tilde{G}_{c,nom}K_r}{(1+D'\tilde{G}_{c,nom}K_r)}$ and $\tilde{T}_2 := \frac{D'\tilde{G}_{c,nom}K_v}{m(1+D'\tilde{G}_{c,nom}K_r)}$. Furthermore, the steady-state tracking error $e_1\triangleq V_\text{ref}-V$ in DC-link voltage is upper bounded by,\\
{\em Centralized case}: $i_\text{ref}=i_\text{load}$
\begin{small}
\begin{equation*}\label{eq:thm2}
	|e_1(j0)|\leq\frac{1}{D'(|K_v(j0)+\eta K_r(j0)|)}|i_\text{ref}(j0)|
\end{equation*}
\end{small}
{\em Decentralized case}: $i_\text{ref}\neq i_\text{load}$
\begin{small}
\begin{eqnarray}\label{eq:thm3}
	|e_1(j0)|\leq &&\frac{|K_r(j0)|}{D'(|K_v(j0)+\eta K_r(j0)|)}|i_\text{ref}(j0)| + \nonumber\\
	&&\frac{D'|K_r(j0)|+1}{D'(|K_v(j0)+\eta K_r(j0)|)}|i_\text{load}(j0)|\nonumber
\end{eqnarray}
\end{small}
\end{theorem}
{\bf Remark 1:} If the steady-state tracking error in DC-link voltage is zero, i.e., $|e_1(j0)| = 0$, then (\ref{eq:thm1}) reduces to the following constraint:
\begin{small}
\begin{equation*}
	\frac{|D_k'i_{L_k}(j0)|}{|D_l'i_{L_l}(j0)|} =  \frac{\gamma_k}{\gamma_l} .
\end{equation*}
\end{small}
i.e., the closed-loop multi-converter system achieves output power sharing given by $|D_1'i_{L_1}(j0)|:\hdots:|D_m'i_{L_m}(j0)| = \gamma_1:\hdots:\gamma_l$. In practice the tracking error $e_1$ is never exactly zero, however, the tracking error is made practically non-existent through an appropriate choice of large DC-gain controllers $K_v, K_r$ resulting from the $\mathcal{H}_\infty$ optimization problem in (\ref{eq:hinf}). Moreover, the design of the controllers is such that $|K_v(j0)|<|K_r(j0)|$ resulting in $|\tilde{T}_1(j0)|\leq1$ and  $|\tilde{T}_2(j0)|\leq1$. The readers are encouraged to refer to the next section for more details where we provide controller transfer functions for few example scenarios.
\\
{\bf Remark 2:} For the {\em decentralized} implementation, it is required that each converter can measure its own inductor current $i_{L_k}$ and DC-link voltage $V$ only.\\
Proof: See Appendix.

\section{CASE STUDIES: SIMULATIONS AND DISCUSSIONS}\label{sec:simulations}
In this section, we report some simulation studies that cover different aspects of the proposed distributed control design. All simulations are performed in MATLAB/Simulink using SimPower/SimElectronics library. Note that the experiments are underway and therefore not reported in this paper. In order to include nonlinearities associated with real-world experiments, and effects of switching frequencies on voltage regulation and power sharing, we use {\em non-ideal} components (such as diodes with non-zero forward-bias voltage, IGBT switches, stray capacitances, parametric uncertainties) and switched level implementation.
\begin{figure*}[!t]
	\begin{center}
	\begin{tabular}{cc}
	\includegraphics[width=0.92\columnwidth]{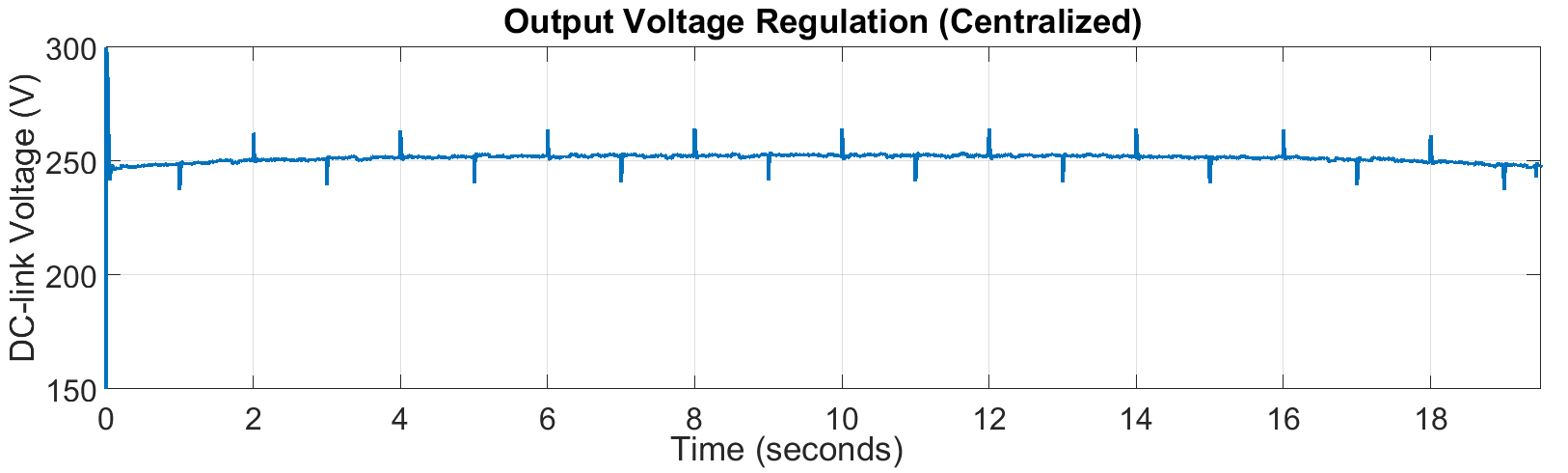}&\includegraphics[width=0.86\columnwidth]{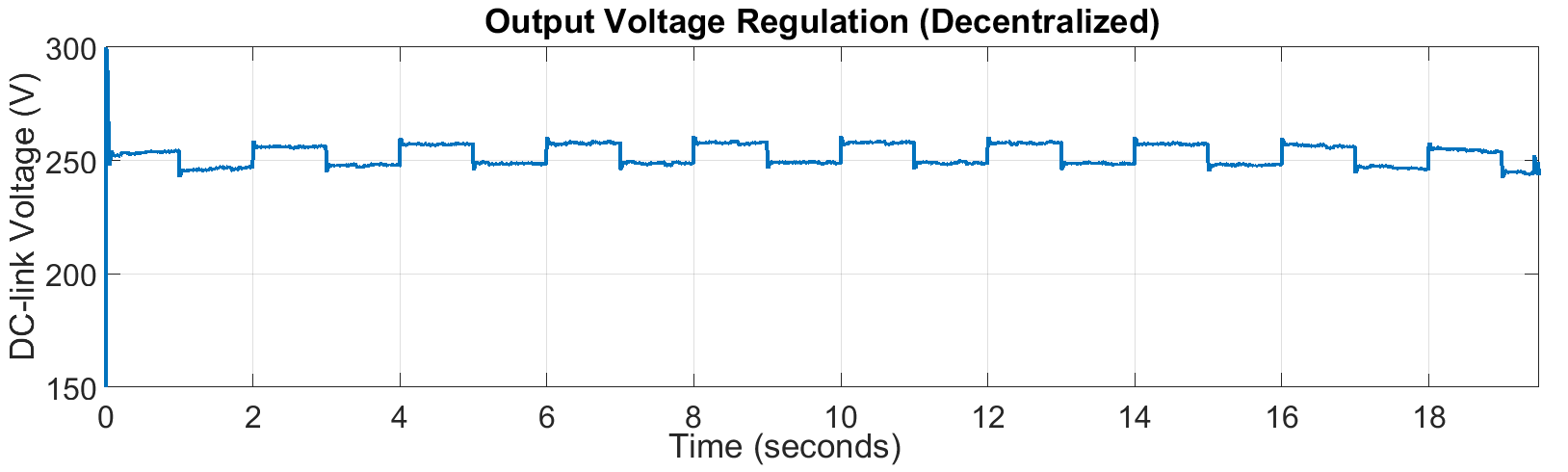}\cr
	(a)&(b)\end{tabular}
	\vspace{-0.5em}
	\caption{{\small Output voltage regulation at the DC-link. (a) Centralized case: Sharp periodic changes in load ($4kW$) result in periodic spikes in output voltage. However, the voltage is soon regulated to $V_\text{ref}=250V$. (b) Decentralized case: In the absence of $i_\text{load}$ measurement, the voltage droops by a prespecified value to accommodate for mismatch in $i_\text{ref}$ and $i_\text{load}$.}}
	\label{fig:voltage_reg}
	\vspace{-1em}
	\end{center}
\end{figure*}
\begin{figure*}[!t]
	\begin{center}
	\begin{tabular}{cc}
	\includegraphics[width=0.92\columnwidth]{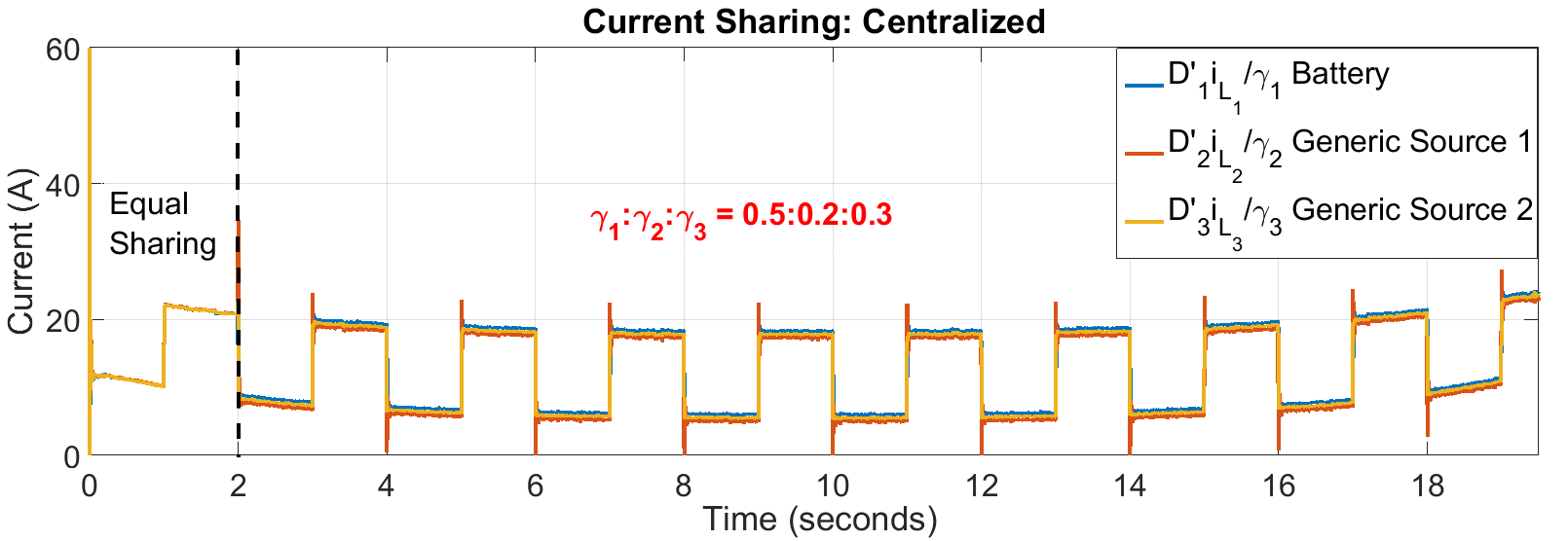}&\includegraphics[width=0.92\columnwidth]{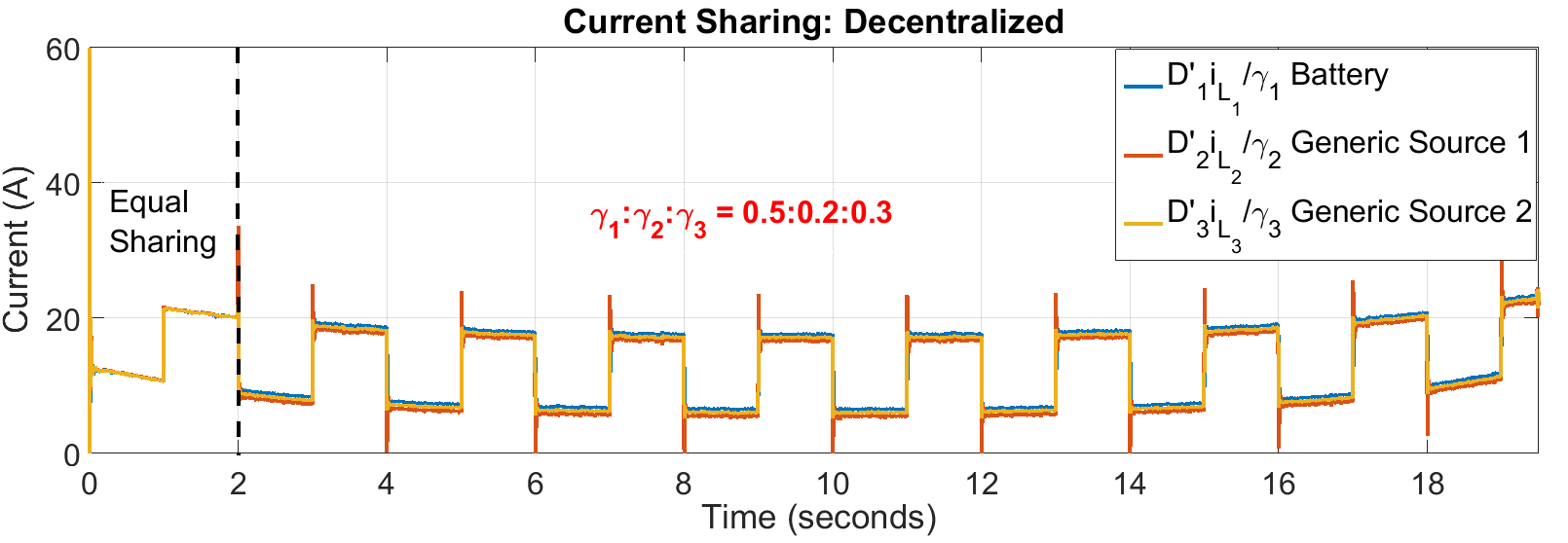}\cr
	(a)&(b)\end{tabular}
	\vspace{-0.5em}
	\caption{{\small Current sharing among three parallel converters. We plot scaled values of the output currents. Overlapping plots of the scaled currents validate that the time-varying sharing requirements are met.}}
	\label{fig:current_sharing}
	\vspace{-1em}
	\end{center}
\end{figure*}
\begin{figure*}[!t]
	\begin{center}
	\begin{tabular}{cc}
	\includegraphics[width=0.92\columnwidth]{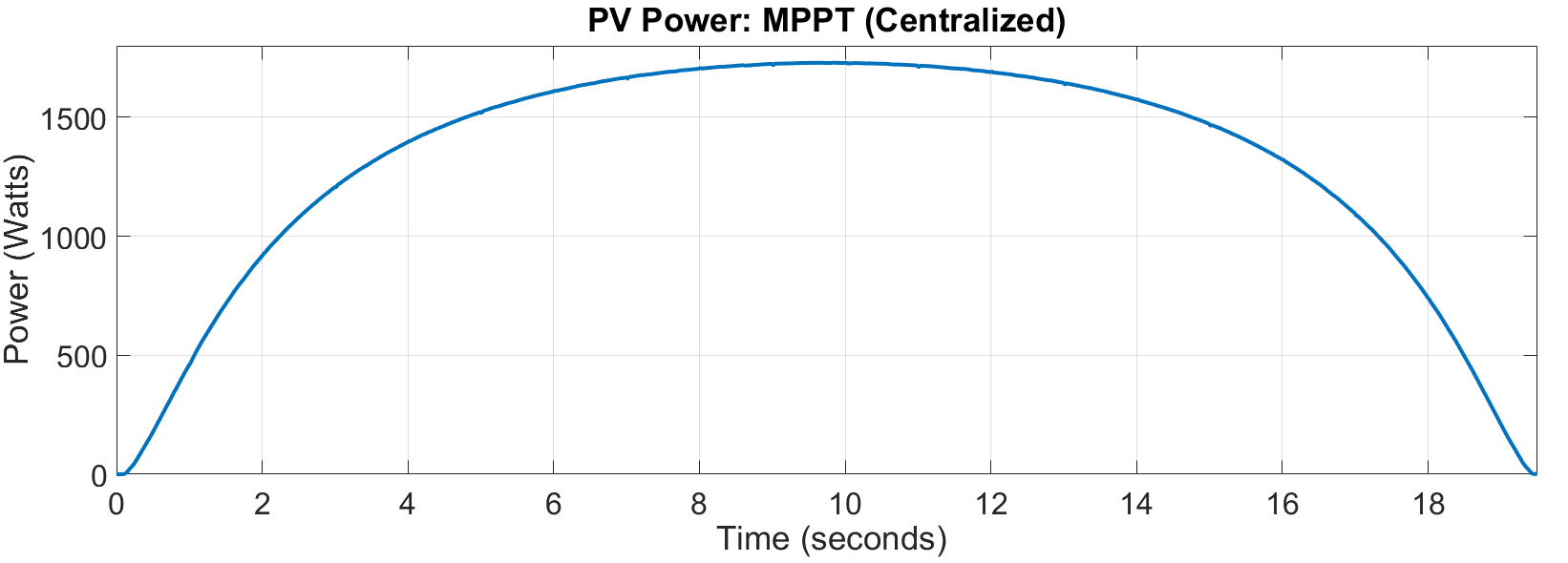}&\includegraphics[width=0.92\columnwidth]{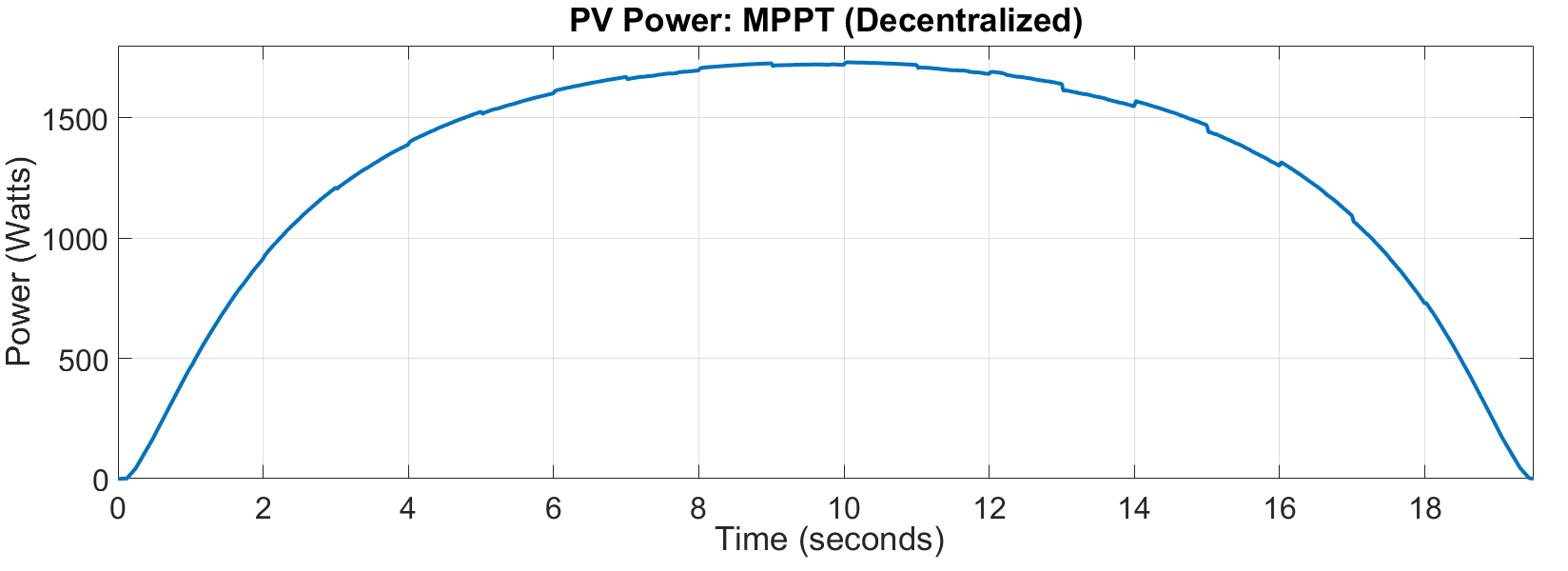}\cr
	(a)&(b)\end{tabular}
	\vspace{-0.5em}
	\caption{{\small PV curve for MPPT algorithm. (a) Centralized case (b) Decentralized case: There are slight periodic jumps in PV curve resulting from voltage droop.}}
	\label{fig:MPPT}
	\vspace{-2em}
	\end{center}
\end{figure*}
For simulations, we consider a parallel network of three boost converters powered by a Li-ion battery and two generic sources (nominal voltages $125V$), respectively. The effectiveness of the proposed control design is well illustrated by considering a challenging practical scenario with ({\em unknown}) fast time-varying load, large uncertainties in inductance and capacitance values, sensor noises and variation in input source voltages. Specifically, we consider the following simulation parameters:
\begin{tcolorbox}
	Source and Converter Parameters:
	\begin{itemize}
		\item Battery State-Of-Charge: $120\% (\sim 135V)$
		\item Generic Sources: $125V$ and $130V$
		\item Inductance: $(L_1, L_2, L_3) = (.096,.12,.14)$mH
		\item Capacitance: $C = 400\mu$F
	\end{itemize}
\end{tcolorbox}
\begin{tcolorbox}
	Photovoltaic System:
	\begin{itemize}
		\item $8$hrs of Insolation data is squeezed to $19.5s$
		\item Peak Power: $1700$W
	\end{itemize}
\end{tcolorbox}
\begin{tcolorbox}
	Output Requirements:
	\begin{itemize}
		\item Loading Conditions: $5kW\pm (2kW@1Hz)$
		\item Power Sharing Requirements:\begin{enumerate}
			\item $(0.33:0.33:0.33), \quad t<2s$
			\item $(0.5:0.2:0.3), \quad 2s\leq t\leq19.5s$
		\end{enumerate}
		\item Desired Output Voltage: $V_\text{ref} = 250V$
	\end{itemize}
\end{tcolorbox}
\begin{tcolorbox}
	Other Simulation Parameters:
	\begin{itemize}
		\item Total Simulation Time: $19.5$s
		\item Droop-Gain Coefficient: $\eta = 1.2667$
		\item Voltage-Sensor Noises:\begin{enumerate}
			\item DC-Offset: $+2V$, Noise: SNR$\sim0.05\%$
			\item DC-Offset: $-2V$
			\item DC-Offset: $+3V$
		\end{enumerate}
	\end{itemize}
\end{tcolorbox}

To illustrate the robustness of the proposed approach, the nominal (or equivalent single converter) inductance, capacitance and steady-state complementary duty-cycle are chosen as $L = 0.12$mH, $C = 500\mu$F and $D' = V_g/V_\text{ref} = 0.5$. The design parameters for the inner-controller $K_c$ are: Damping factors $\zeta_1 = 0.7$, $\zeta_2 = 2.2$, and bandwidth $\tilde{\omega} = 2\pi 300$rad/s.

The outer controllers $K_{v}$ and $K_{r}$ are obtained by solving the stacked $\mathcal{H}_\infty$ optimization problem (see Eq. (\ref{eq:hinf}))\cite{skogestad2007multivariable} with the weighting functions:
\begin{small}
	\begin{eqnarray}
		W_1 = \frac{0.4167(s+452.4)}{(s+1.885)} \quad W_2 = \frac{0.4167(s+1056)}{(s+4.398)}\nonumber \\
		W_3 = 0.04 \qquad W_4 = \frac{37.037(s+314.2)}{(s+3.142\times10^4)}\nonumber
	\end{eqnarray}
\end{small}
The resulting outer-controllers are reduced to sixth-order using balanced reduction \cite{dullerud2013course} and are given by,
\begin{small}
	\begin{eqnarray}
		K_{v} &=& \frac{-0.00064(s-4.615e9)(s+6007)}{(s+1.604e4)(s+578.3)}\nonumber \\
		&&\frac{(s^2+5042s+5.97e6)}{(s^2+1061s+5.69e5)}\frac{(s^2+753.6s+1.039e5)}{(s^2+7.354e4s+2.074e9)}\nonumber\\
		K_{r} &=& \frac{0.00267(s+181.3)(s+0.001012)}{(s+4.395)(s+0.001013)}\nonumber \\
		&&\frac{(s^2+5141s+6.065e6)}{(s^2+1059s+5.694e5)}\frac{(s^2+3.818e6s+2.804e11)}{(s^2+9.783e4s+3.69e9)} \nonumber 
	\end{eqnarray}
\end{small}
{\em Inclusion and characterization of PV module}: Photovoltaics are technically treated as {\em current} sources. In a microgrid setup, a PV module is interfaced with the DC-link through a Boost converter and is controlled using the maximum power point tracking (MPPT) algorithm. The output current of PV $i_\text{PV}$ is directly proportional to the (time-varying) irradiance and is included in our proposed formulation by regarding $i_\text{PV}$ as part of the disturbance signal, i.e., the net disturbance current is modeled as $i_\text{load}-i_\text{PV}$. In this simulation study, we squeeze worth $8$ hours of insolation data into a total duration of $19.5s$ amounting to rapidly varying irradiance (and hence the disturbance current $i_\text{PV}$). Note that this is a challenging problem to consider from not only MPPT point-of-view, but also on the account of DC-link voltage regulation coupled with time-varying (and uncertain) loads and time-dependent power-sharing requirements.

{\em Results}: The controllers derived for the nominal single converter system are then extended for a parallel multi-converter design as described in Sec. \ref{sec:Control_Many}. The initial DC-link voltage is considered at $0V$. Fig. \ref{fig:voltage_reg} shows the voltage regulation at the DC-link to the reference $V_\text{ref}=250V$ for the centralized ($i_\text{load}$ measurement available) and decentralized implementations. Note that the DC-link load changes by $4kW$ every second ($3kW$ to $7kW$, and $7kW$ to $3kW$). The reference current is considered as $i_\text{ref}=5kW/250V = 20A$. While in the centralized case (Fig. \ref{fig:voltage_reg}a), the voltage is maintained at $V_\text{ref}=250V$ with small periodic spikes attributed to sudden load changes, the decentralized implementation results in controlled voltage droop of $10V$ peak-to-peak around the desired DC-link voltage. 

Fig. \ref{fig:current_sharing} presents the results for time-varying sharing. For ease of illustration, the scaled output currents $D'i_L/\gamma$ are plotted. The power sharing requirements are described in the output requirements section. Overlapping values of scaled currents depict excellent sharing performance; moreover, the control design achieves step-change in power sharing seamlessly at $t=2s$.

Fig. \ref{fig:MPPT} shows the result of incremental conductance based MPPT for the chosen PV system for the two scenarios. It is observed that the availability of $i_\text{load}$ measurement results in smoother PV curve in Fig. \ref{fig:MPPT}a, while the decentralized implementation has slight periodic jumps in PV curve resulting from voltage droop.

\section{CONCLUSIONS AND FUTURE WORK}\label{sec:conclusion}
In this work, we propose a distributed control architecture for voltage tracking and power sharing for a network of DC-DC converters connected in parallel. The proposed design is capable of achieving multiple objectives such as robustness to modeling uncertainties, reference DC voltage generation and output power sharing among multiple DC sources. The controllers are designed using a robust optimal control framework. A remarkable feature of the proposed control approach is the ability to achieve DC-link voltage regulation for both centralized and decentralized scenarios without the need to modify the controller structure or the parameters. Experimental validation of the proposed approach is underway and will soon be reported in our subsequent work.

\section*{ACKNOWLEDGMENT}
The authors would like to acknowledge ARPA-E NODES program for supporting this work.

\section*{APPENDIX}
\subsection*{Sketch of Proof of Theorem 1: System Equivalence}
\begin{proof}
	The {\em equivalence} is a direct consequence of cleverly chosen architecture. Note that for the single converter system in Fig. \ref{fig:single_block} with $\tilde{G}_c(s) = \tilde{G}_{c,nom}(s)$, the error signal $e_2$ (input to controller $K_r$) is given by
	\begin{small}
		\begin{eqnarray}\label{eq:proof1_eq1}
			e_2 &=& [i_\text{ref}+\eta e_1-D'\tilde{G}_{c,nom}\hat{u}]\nonumber \\
			\Rightarrow e_2 &=& i_\text{ref}+(\eta-D'\tilde{G}_{c,nom})e_1-D'\tilde{G}_{c,nom}K_re_2.
		\end{eqnarray}
	\end{small}
	For the multi-converter system in Fig. \ref{fig:model_many}b, we have $e_1^{(k)}=V_\text{ref}-V:=e_{1}$. Let us denote the total error in current mismatch $\sum_{k=1}^me_2^{k}$ by $e_{2}$. Therefore, from Fig. \ref{fig:model_many}b,
	\begin{small}
		\begin{eqnarray}\label{eq:proof1_eq2}
			e_2^{(k)} = \gamma_k[i_\text{ref}+\eta e_1]-D'\tilde{G}_{c,nom}\left(\frac{1}{m}K_{v}e_1+K_re_2^{(k)}\right)\nonumber \\
			\underbrace{\sum_{k=1}^me_2^{(k)}}_{e_{2}} = \sum_{k=1}^m\gamma_k[i_\text{ref}+\eta e_1] - D'\tilde{G}_{c,nom}\Big(K_ve_1 - K_r\underbrace{\sum_{k=1}^me_2^{(k)}}_{e_{2}}\Big).
		\end{eqnarray}
	\end{small}
	Using the fact that $\sum_{k=1}^m\gamma_k=1$, the above equation reduces to (\ref{eq:proof1_eq1}). Similarly, the expression for tracking error in voltage $V_\text{ref}-V$ is identical for the single and multiple converters case. Moreover for the multi-converter system, the output voltage at the DC-link is given by $V = G_v(-i_\text{load}+D'\tilde{G}_{c,nom}(K_ve_{1}+K_re_{2}))$. Since the expressions for $e_1$ and $e_2$ are identical for the single and multiple converters case and are written in terms of the exogenous variables $V_\text{ref}, i_\text{ref}, i_\text{load}$, the corresponding transfer functions from the exogenous variables to the DC-link voltage $V$ are also identical, and hence establishes the required {\em equivalence}.
	Similar conclusions can be drawn for other signals, such as DC-link voltage $V$, and hence establishes the required {\em equivalence}.
\end{proof}
	
\subsection*{Sketch of Proof of Theorem 1: Power Sharing}
\begin{proof}
	From (\ref{eq:proof1_eq2}), the error in current reference for the $k^{th}$-converter is given by
	\begin{small}
		\begin{eqnarray}\label{eq:proof2_eq1}
			e_2^{(k)} = \left(\frac{\gamma_k}{1+D'\tilde{G}_{c,nom}K_r}\right)i_\text{ref} + \left(\frac{\gamma_k\eta-\frac{D'}{m}\tilde{G}_{c,nom}K_v}{1+D'\tilde{G}_{c,nom}K_r}\right)e_{1}.
		\end{eqnarray}
	\end{small}
	From Fig. \ref{fig:model_many}b, the output current $i_{k}=D_k'i_{L_k}$ of the $k^{th}$ converter is given by
	\begin{small}
		\begin{eqnarray}\label{eq:proof2_eq2}
			i_k = D'\tilde{G}_{c,nom}\left[\frac{1}{m}K_ve_{1}+K_re_2^{(k)}\right].
		\end{eqnarray}
	\end{small}
	Thus from (\ref{eq:proof2_eq1}) and (\ref{eq:proof2_eq2}), we obtain
	\begin{small}
		\begin{equation*}
		\resizebox{.5 \textwidth}{!} 
		{
			$i_k = D'\tilde{G}_{c,nom}\left[\left(\frac{\gamma_kK_r}{1+D'\tilde{G}_{c,nom}K_r}\right)i_\text{ref}+\left(\frac{\frac{1}{m}K_v+\eta\gamma_kK_r}{1+D'\tilde{G}_{c,nom}K_r}\right)e_1\right]$.
		}
		\end{equation*}
	\end{small}
	Therefore, we have
	\begin{small}
		\begin{equation*}
			\left|\frac{i_k(j0)}{\gamma_k}-\frac{i_l(j0)}{\gamma_l}\right|\leq \left(\eta|\tilde{T}_1(j0)| + \left|\frac{1}{\gamma_k}-\frac{1}{\gamma_l}\right||\tilde{T}_2(j0)|\right)|e_1(j0)|
		\end{equation*}
	\end{small}
	The expressions for the bounds on the tracking error for the two scenarios is directly obtained from (\ref{eq:voltage_single}) and the {\em system equivalence} described earlier.
\end{proof}

\bibliographystyle{IEEEtran} 
\bibliography{myRef}

\end{document}